\documentclass{amsart}
\usepackage{graphicx}
\usepackage{color}
\theoremstyle{definition}

\theoremstyle{remark}

\numberwithin{equation}{section}

\newcommand{\regT}{{\rm T}}
\newcommand{\regF}{{\rm F}}
\begin{document}

\title[A Symbolic Finite-state approach ]{A Symbolic Finite-state
approach for Automated Proving
of Theorems in Combinatorial Game Theory}%
\author{Thotsaporn  ``Aek'' Thanatipanonda}%
\author{Doron Zeilberger}%
\address{Department of Mathematics}%
\address{Hill Center-Busch Campus}%
\address{Rutgers University}%
\address{110 Frelinghuysen Rd}%
\address{Piscataway, NJ 08854-8019}%
\address{USA}%

\email{thot@math.rutgers.edu,zeilberg@math.rutgers.edu}%

\subjclass{91A46}%
\keywords{Combinatorial Game Theory}%


\begin{abstract}
   We develop a finite-state automata
approach, implemented in a Maple package
{\tt ToadsAndFrogs} available from our websites,
for conjecturing, and then
rigorously proving, values for
large families of positions in Richard
Guy's combinatorial game ``Toads and Frogs''.
In particular, we prove a conjecture of Jeff Erickson.
\end{abstract}

\maketitle
\section{Introduction}

The game {\it Toads and Frogs}, invented by Richard Guy,
is extensively discussed in ``Winning Ways''[1], the famous
classic by Elwyn Berlekamp, John Conway, and Richard Guy, that is
the {\it bible} of combinatorial game theory. \\

This game got so much coverage because of the
simplicity and elegance of its rules, the beauty of its
analysis, and  as an example of a  combinatorial game
whose positions do not always have values that are numbers. \\

The game is played on a $1 \times n$ strip with either
Toad(T) , Frog(F)
or $\Box$ on the squares. Left plays T and Right plays F.
T may move to the immediate square on its right, if it happens
to be empty,
and F moves to the next empty square on the
left, if it is empty.
If T and F are next to each other, they have an option to jump over one
another,
in their designated directions, provided they lend on an empty square.
(See [1], page 14).\\

\noindent In symbols: the following moves are legal for Toad: \\

\noindent $\dots \regT \Box \dots \;\ \rightarrow \;\ \dots \Box \regT
\dots
$,\\
$\dots \regT \regF \Box \dots \;\ \rightarrow \;\ \dots \Box \regF \regT
\dots $
\;\ , \\

\noindent and  the following moves are legal for Frog: \\

\noindent $\dots \Box \regF \dots \;\ \rightarrow \;\ \dots \regF \Box
\dots
$,\\
$\dots \Box \regT \regF \dots \;\ \rightarrow \;\ \dots \regF \regT \Box
\dots
\quad $ . \\

Already in ``Winning Ways''[1], there is
some analysis of Toads and Frogs positions, but on {\it specific}, small
boards, such as $\regT \regT \regT \Box \regF \regF$. In 1996,
Jeff Erickson[2]
analyzed more general positions. At the end he made five
conjectures about the values of some
families of positions. All of them are ``starting''
positions
(i.e. positions where all T's are rightmost and all F's are leftmost).\\

To be able to understand the present article,
readers need some knowledge of
combinatorial game theory, that can be found in [1].
In particular, readers should be familiar with the
notion of {\it value} of a game. Recall that values
are not always numbers (not even surreal ones).\\

Let's recall the {\it bypass
reversible move rule},
the {\it dominated options rule} (see [1] page 62-64) and Erickson's
{\it Terminal Toads Theorem} (see[2]).\\

\noindent \textbf{Bypassing right's reversible move rule}\\

\begin{figure} [ht]
\input{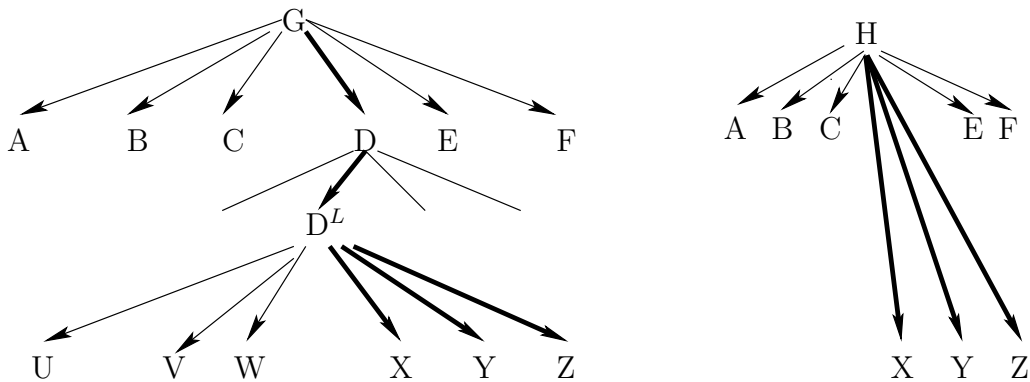}
\caption{The Bypassing reversible move rule.}
\end{figure}

\noindent $G = H$ if $D^{L} \geq G.$\\

\noindent \textbf{The Dominated options Rule} \\

\noindent Let $G = \{A,B,C,... \mid D,E,F,.. \}.$\\
If $A \geq B$ and $D \geq E$ then $G = \{A,C,... \mid E,F,.. \}.$\\

\noindent \textbf{The Terminal Toads Theorem}: Let X be any position.
Then\\
$ {\rm X} \regT \Box^{n} = {\rm X} \Box^{n} + n.$\\

The only notation we use is $ * \;\ ( = \{ 0 \mid 0 \})$. We will not use
any shorthand notation like $\uparrow$, $\Uparrow$, etc.\\

Next, we will explain the
method through examples, and describe how to implement the method
when applied to certain classes of positions. Finally, we discuss a new
conjecture and possible future work. \\

Everything is fully implemented in a  Maple package,
{\tt ToadsAndFrogs}, written by the first author,
available from our websites.

\section{A Symbolic Finite-State Method}

\noindent We define two classes of positions:\\

\noindent Class {\bf A}: All the positions that have a {\it fixed}
number of occurrences of $\Box$ and F,
but a {\it variable} (symbolic) number of T's in-between the $\Box$'s and
F's.
\\

\noindent class {\bf B}:
All the positions that have a {\it fixed}
number of occurrences of T's and F's,
but a {\it variable} (symbolic) number of $\Box$'s in-between the T's and
F's.
\\

\noindent $Aij$ := the class in which we have exactly $i$ occurrences of
$\Box$
and exactly $j$
occurrences of $\regF$.\\

\noindent $Bij$ :=
the class in which we have exactly $i$ occurrences of $\regT$ and exactly
$j$
occurrences
of $\regF$.\\

For any {\it specific} position,
we can always
compute its value,  by using the recursive definition of the value.
But this is mere {\it number-crunching}.
After collecting enough data, and examining it,
if we are lucky,
we (or rather our computers) can detect a uniform {\it pattern}, and {\bf
conjecture}
an {\bf explicit} formula for the values of the studied family,
in terms of the symbolic parameters.
Once conjectured, these conjectured explicit expressions
can be proved  by induction on the symbolic
parameters. The beauty and novelty of
our approach is that everything is done {\bf automatically}.
First the {\it conjecturing} parts,
but more dramatically, the {\it proving} part.
We teach the computer how to conjecture, by looking for general
patterns, and then how to use induction in order to prove its own
conjectures.
\\

This {\it activity} of {\bf computer-generated} mathematics is in sharp
contrast to the traditional approach of [2], that merely uses
the computer as a calculator, to generate numerical data,
and everything else, the conjecturing, and the
proving (when feasible) is done by humans. \\

We believe that the present methodology is of potential use
in many other branches of mathematics, and ``Toads and Frogs''
is but an instructive arena for presenting a general approach for
computer-generated research.\\

When we analyze each class of positions, we are naturally lead,
by the recursive definition of the {\it value} (of a game), to other
classes of positions. Luckily, at least in all the
cases encountered so far,
there are always a {\bf finite}
number of different classes, that we can name
``symbolic states''. If the
(symbolic) value of each ``state'' in the class is conjectured to have
a (symbolic) explicit expression,
then we can prove the truth of {\it all} these conjectures
{\bf all at once} by applying
induction on the recurrence relations.
Note that in order for this to work we need to conjecture explicit
expressions for {\bf all} the states, so we usually get much more
than we bargained for.
\\

We will demonstrate the method  with the two simplest nontrivial
classes: A11 and B11.\\

\noindent \textbf{First example}: Type A11: one $\Box$ and one F \\

\noindent Let $f(a,b)$ be the value of $\regT^{a} \Box \regT^{b} \regF$.
\\
Let $g(a)$ be the value of $\regT^{a} \regF \Box$.  \\

Here, of course, $\regT^{a}$ means $\regT$ repeated $a$ times, so the
`game'
$f(a,b)$, for example, stands for a doubly-infinite set of
starting positions.\\

\noindent \textbf{Recurrences}: \\
Note that if any parameter of the function is negative then it return
NULL.\\

\begin{figure} [ht]
\input{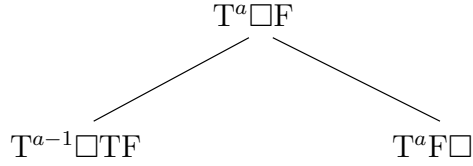}
\caption{Recurrence for $f(a,0), a \geq 0$.}

\bigskip
$f(a,0) = \{ f(a-1,1) \mid g(a)\} ,  a \geq 0. $\\

\end{figure}

\begin{figure} [ht]
\input{tf12.pstex_t}
\caption{Recurrence for $f(a,1), a \geq 0$.}

\bigskip
$f(a,1) = \{ f(a-1,2) \mid g(a)+1\} ,  a \geq 0.$ \\
\end{figure}

\begin{figure} [ht]
\input{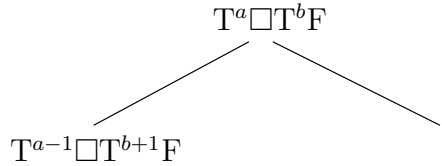}
\caption{Recurrence for $ f(a,b), a \geq 0, b \geq 2 $.}

\bigskip
$f(a,b) = \{ f(a-1,b+1) \mid \;\ \} ,  a \geq 0 , b \geq 2.$ \\
\end{figure}

\newpage

\begin{figure} [ht]
\input{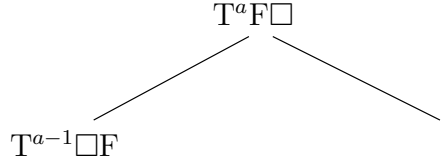}
\caption{Recurrence for $g(a), a \geq 0$.}

\bigskip
$g(a) = \{ f(a-1,0) \mid \;\ \},  a \geq 0$. \\
\end{figure}

The above recurrences can be easily used to
crank out {\it numerical data} for
small (and not so small) values of $a$ and $b$.
Then the computer
{\it automatically} makes the following {\it symbolic}
conjectures.\\

\noindent \textbf{Conjectures}:
\[ \begin{array}{llll}
f(0,0) &=& -1. \\
f(a,0) &=& \{\{a-2 \mid 1 \} \mid 0 \} &, \;\ a \geq 1. \\
f(a,1) &=& \{ a-1 \mid 1 \} &, \;\ a \geq 0.\\
f(a,b) &=& a &,\;\ a \geq 0,  \;\ b \geq 2.\\
g(a) &=& 0 &, \;\ a \geq 0.\\ \\
\end{array} \]

Once conjectured, the proof is routine, and also can (and was!)
done by computer. One checks the obvious initial conditions
and verifies that the above expressions satisfy the above
{\it defining} relations. Indeed, the computer easily
verifies that

\[ \begin{array}{llll}
f(a,0) &=& \{f(a-1,1) \mid g(a) \} = \{\{ a-2\mid 1 \} \mid 0 \} &,\;\  a
\geq 1. \\
f(0,1) &=& \{ \;\ \mid g(0)+1 \} = \{ \;\ \mid 1 \} = 0 = \{ -1 \mid 1 \}.
\\
f(a,1) &=& \{f(a-1,2) \mid g(a)+1 \} = \{ a-1\mid 1 \} &, \;\ a \geq 1. \\
f(0,b) &=& \{ \;\ \mid \;\ \} = 0 &, \;\ b \geq 2. \\
f(a,b) &=& \{f(a-1,b+1) \mid \;\ \} = \{ a-1\mid  \;\ \} = a  &, \;\ a
\geq 1,
b \geq 2. \\
g(0) &=&  \{ \;\ \mid \;\ \} = 0.  \\
g(1) &=& \{f(0,0) \mid \;\ \} = \{ -1\mid \;\ \} = 0. \\
g(a) &=& \{f(a-1,0) \mid \;\ \} = \{\{\{ a-3\mid 1 \} \mid 0
\}\mid \;\ \}\\
     &=& \{ \;\ \mid \;\ \}$ (!! by bypass reversible move rule) $= 0 &,
\;\ a
\geq 2. \\ \\
\end{array} \]

Note that the above values for $f(a,0)$ ($a \geq 1$) agree with the case
$b=1$ of Theorem 5.2 of [2]. \\

\noindent \textbf{Second Example}: Type B11: one T and one F.     \\

\noindent Let $f(a,b,c) := \Box^{a} \regT \Box^{b} \regF \Box^{c}.$ \\

\noindent Now we have a {\it three-} parameter family!\\

\noindent \textbf{Initial Conditions and Recurrences}: \\
\[ \begin{array}{llll}
f(0,0,0) &=& \{ \;\ \mid \;\ \}.\\
f(a,0,0) &=& \{ \;\ \mid (-a+1)+1\} = \{ \;\ \mid -a+2 \} &, \;\ a \geq
1.\\
f(0,0,c) &=& \{ (c-1)-1\mid \;\ \} = \{ c-2 \mid \;\ \} &, \;\ c \geq 1.\\
f(a,0,c) &=& \{ c-a-2\mid c-a+2\} &, \;\ a \geq 1, c \geq 1.\\
\end{array} \]

\begin{figure} [ht]
\input{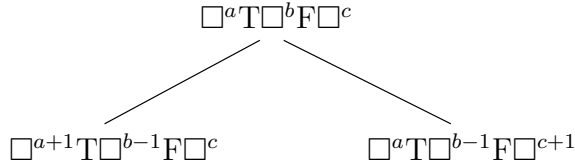}
\caption{Recurrence of $f(a,b,c), a \geq 0, c \geq 0, b \geq 1 $.}

\bigskip
$f(a,b,c) = \{ f(a+1,b-1,c)\mid f(a,b-1,c+1) \}, \;\ a \geq 0, c
\geq 0, b \geq 1$.\\

\end{figure}

By using these recurrences {\it numerically}, the computer
cranks out enough data, that enables it to make the following\\

\noindent \textbf{Conjecture}: \\
\[ \begin{array}{llll}
f(a,b,c) &=& \{ c-a-2\mid c-a+2\} &, \;\ a \geq 0, c \geq 0 , b\mbox{ is
even}.\\
f(a,b,c) &=& \{\{ c-a-3\mid c-a+1\} \mid \{ c-a-1\mid c-a+3\}\} &, \;\ a
\geq 0, c \geq 0, \, b \mbox{ is odd }.\\
\end{array} \] \\

\noindent \textbf{Proof:} by induction: on $b$. \\

\noindent \textbf{Base case}: $b$ = 0\\
We have
\[ \begin{array}{llll}
f(0,0,0) &=& 0 = \{-2 \mid 2 \} .  \\
f(a,0,0) &=&  \{\;\; \mid -a+2 \} = \{-a-2 \mid -a+2 \} , \;\ &a \geq 1.\\
f(0,0,c) &=&  \{ c-2 \mid  \;\ \} = \{c-2 \mid c+2 \} ,&c \geq 1.\\
f(a,0,c) &=& \{ c-a-2\mid c-a+2\} ,&a \geq 1, c \geq 1.\\
\end{array} \]

\noindent \textbf{Induction step} on $b$ \\

\noindent Case 1) $b$ is even and  $b \neq 0$:
\[ \begin{array}{llll}
f(a,b,c) &=& \{ f(a+1,b-1,c)\mid f(a,b-1,c+1) \} ,\;\ a \geq 0, c \geq
0.\\
  &=& \{\{ c-a-4\mid c-a\} \mid \{ c-a-2\mid c-a+2\}\} \\
  && \mid \{\{ c-a-2\mid c-a+2\} \mid \{ c-a\mid
c-a+4\}\}\}.\\
  &=& \{ c-a-2 \mid c-a+2 \}. \\
\end{array} \]

\noindent Case 2) $b$ is odd
\[ \begin{array}{llll}
f(a,b,c) &=& \{ f(a+1,b-1,c)\mid f(a,b-1,c+1) \} &, a \geq 0, c \geq 0.\\
         &=& \{\{ c-a-3\mid c-a+1\} \mid \{ c-a-1 \mid c-a+3\}\}. \\
\\
\end{array} \]

The second example is related to the results of
Erickson [2] as follows.
The case $b=0$ is Lemma 4.1 of [2], while the
case $a=0,c=0$ coincides with  the case $a=1$ of Theorem 5.2.
Note that we need the extra elbow-room of a three-parameter
family to enable the inductive argument. \\

\section{How far can the symbolic finite state method go?}

As we mentioned in the previous section, the finite state method works
perfectly well when the value
of every position in the class has a discernible pattern. This seems
to be the case for class A. We wrote a computer program in Maple to
first {\it calculate}, then {\it conjecture}, and
finally {\it prove}, the values of general positions in class A
automatically. The program now works for positions with any fixed
number of $\Box$'s and with one Frog. For the class where we have more
than one Frog, it is harder to find conjectures, for humans, and
{\it even} for computers. We conjectured some classes with two
Frogs(A12, A22, A32)  by
hand and put it in the computer program to prove the conjectures.\\

The list of the results for the classes A11, A21, A31, A41, A51,  A12,
A22, A32 can be found in both authors' websites. \\

As a very special case of our results for the class A32, we
get a proof of
Erickson's[2] conjecture 2, that claims that the value of
$\regT^{a} \Box\Box\Box \regF \regF$  is $ \{ a-2 \mid a-2 \}$,  ( $a \geq
2$).\\

In [3], the first-named author
of the present paper discusses the value of {\it any} position with one
$\Box$ and {\it any} number of Toads  and Frogs
(Therefore we are done with class $A1n$,
$n \geq 1$). This general class with one $\Box$ is the only general class
we are able to figure out the patterns for.\\

\noindent We now turn our attention to class B. We solved class B11 in the
previous section. For B21: TTF, we already have a difficulty. The formulas
in this class are long and hard to find in a canonical form. We will
discuss
this in the appendix.\\

\section{A Conjecture and Future Work}

\noindent \textbf{Conjecture}: \\

1) We always have ``nice compact'' formulas for every position in class
A.\\

\noindent \textbf{Future Work}: \\

1) Implement the symbolic finite state method for the class B21.\\

2) We have seen systems of recurrence relations arising naturally
in each class. We solved the recurrences by ``guessing''
(automatically, of course) the answers
(using predefined {\it ansatzes}) and then proving
them by induction. It would be interesting to develop general
algorithms for systematically solving the recurrnces, without
the need for ``guessing''. \\

\renewcommand{\theequation}{A-\arabic{equation}}
\setcounter{equation}{0}  

\section*{APPENDICES}  
\appendix{}

\section{On the difficulty of class B21: TTF}
\noindent \textbf{B21}: TTF

\noindent $f(a,b,c,d) := \Box^{a} \regT \Box^{b} \regT \Box^{c} \regF
\Box^{d}.$\\
$g(a,b,c) := \Box^{a} \regT \Box^{b} \regF \Box^{c}.$\\

\noindent We already knew the solution of g since it is exactly B11.\\
We can now focus on f.\\

\noindent \textbf{Recurrences}:
\[ \begin{array}{llll}
f(a,0,0,0) &=& \{ \;\ \mid \;\ \} = 0.\\
f(a,0,0,d) &=& \{ g(a,1,d-1)+d-1 \mid \;\ \} \\
&=& \{\{\{ d-a-4 \mid d-a \} \mid \{ d-a-2 \mid d-a+2 \}\}
\mid \;\ \} \\
&&,\;\ a\geq 0, \;\ d \geq 1.\\
f(a,b,0,0) &=& \{ f(a+1,b-1,0,0) \mid g(a,b-1,1)+1 \} \\
&&, \;\ a \geq 0, b \geq 1.\\
f(a,b,0,d) &=& \{ f(a+1,b-1,0,d) ,\;\ g(a,b+1,d-1)+d-1  \mid
g(a,b-1,d+1)+d+1
\} \\
&&, \;\ a \geq 0, b \geq 1, d \geq 1.\\
f(a,b,c,d) &=& \{ f(a+1,b-1,c,d) ,\;\ f(a,b+1,c-1,d)  \mid f(a,b,c-1,d+1)
\}
\\
&&, \;\ a \geq 0, b \geq 0, c \geq 1, d \geq 0.\\
\end{array} \]

\noindent Note: $f(a,0,0,d)$ has been discussed before as lemma 4.3 by
Erickson.\\

\noindent \textbf{A nice formula for $f(a,b,0,0)$}. \\

\noindent For b =1:
\[ \begin{array}{llll}
f(a,1,0,0) &=& \{ f(a+1,0,0,0) \mid g(a,0,1) +1 \} \\
&=& \{ 0 \mid \{ -1-a \mid 3-a \} +1 \} \\
&=& \left\{ \begin{array}{ll}
   \frac{1}{2}  & \mbox{, a = 0,1} \\
   \{0 \mid 3-a \} & \mbox{, a} \geq \mbox{2} \\
   \end{array}
   \right.
\end{array} \]

\noindent For b $\geq$ 2 and b is even:
\[ \begin{array}{llll}
\noindent f(a,b,0,0) &=& \left\{ \begin{array}{ll}
   1  & \mbox{, a = 0 } \\
   -a+2 & \mbox{, a } \geq 1. \\
   \end{array}
   \right. \\
&=& \{ \;\ \mid a \} -a+2. \\
\end{array} \]

\noindent For b $\geq$ 2 and b is odd. \\
\[ \begin{array}{llll}
f(a,b,0,0) = \left\{ \begin{array}{ll}
   \{1 \mid 1\}  &, a = 0 \\
   \frac{1}{2} &,  a = 1 \\
   -a+2  &, a \geq 2 .\\
   \end{array}
   \right.
\end{array} \] \\

However for $f(a,b,0,d) , a\geq 0, b \geq 1, d\geq 1$ , the formulas get
longer and longer and we started to lose track of them, and consequently
failed to find formulas in this case. It should be possible to write
Maple code specifically to find a pattern for the values of positions in
class B. The authors expect the formulas in other classes of type B (for
example B22: TTFF) to be even more complicated than B21, since it has to
build
up from B21. \\

It appears that the positions in class B have periodicity and they
need more care to formulate the right conjectures. \\

\section{ About the program}

Our Symbolic Finite-State Method was implemented by one of us (TT)
in Maple.
He first wrote a program to recursively calculate the values of games.
Then he improved the program by making use of the symbolic computation
capability of Maple, to formulate
conjectures, and prove the values of game-positions.
The whole proof process was completely automated. Below is the short
description of the program. See the authors' web sites for complete
details of the program.\\

\noindent \textbf{ToFr} \\

{\bf Input}: the specific position of the game. \\

{\bf Output}: the value of the game in canonical form.\\

\noindent \textbf{SVG} \\

{\bf Input}: the value of the game, could be symbolic. \\

{\bf Output}: the value of the game in canonical form.\\

Note: This program can also be used for other combinatorial games.\\

\noindent \textbf{MainConj} \\

{\bf Input}: number of $\Box$ and number of F. \\

{\bf Output}: The list of conjectures. \\

\noindent \textbf{Prove} \\

{\bf Input}: number of $\Box$ and number of F.\\

{\bf Output}: the values of all of the positions in this class.\\

The program currently only works for one Frog with any fixed number of
$\Box$. With more than one Frog, it gets harder to find
conjectures. But one could find conjectures by hand and feed them to the
subfunctions in Prove. The program can help verify such
humanly-made conjectures.\\

\noindent Obviously, there is still a lot of work to be done, but
let's remember that \\

\noindent ``
{\bf Every great artwork always starts from a rough draft}''.\\


\end{document}